\documentclass[12pt]{article}
\usepackage{amsmath,amsfonts,amsthm,amssymb,eucal}
\newcommand{\Qp}{\mathbb Q_p}
\newcommand{\DA}{D^\alpha}
\newcommand{\ufrac}{\frac{1-p^{-\alpha}}{1-p^{\alpha-1}}}
\newcommand{\iafi}{I^{\alpha} f }
\begin{document}
\newtheorem{lem}{Lemma}
\newtheorem{teo}{Theorem}
\newtheorem{prop}{Proposition}
\newtheorem*{defin}{Definition}
\pagestyle{plain}
\title{Asymptotic Properties of the $p$-Adic Fractional Integration Operator}
\author{ \textbf{Anatoly N. Kochubei}\\
\footnotesize Institute of Mathematics,\\
\footnotesize National Academy of Sciences of Ukraine,\\
\footnotesize Tereshchenkivska 3, Kyiv, 01004 Ukraine,\\
\footnotesize E-mail: kochubei@imath.kiev.ua
\and
\textbf{Daniel S. Soskin}\\
\footnotesize Faculty of Mathematics and Mechanics,\\
\footnotesize Taras Shevchenko Kyiv National University,\\
\footnotesize Volodymyrska 64, Kyiv, 01033 Ukraine,\\
\footnotesize E-mail: dssoskin@gmail.com}

\date{{\it To the blessed memory of M. L. Gorbachuk}}

\maketitle

\bigskip
\begin{abstract}
We study asymptotic properties of the $p$-adic version of a fractional integration operator introduced in the paper by A. N. Kochubei, Radial solutions of non-Archimedean pseudo-differential equations, {\it Pacif. J. Math.} {\bf 269} (2014), 355--369.
\end{abstract}
\vspace{2cm}
{\bf Key words: }\ $p$-adic numbers; Vladimirov's $p$-adic fractional differentiation operator; $p$-adic fractional integration operator; asymptotic expansion

\medskip
{\bf MSC 2010}. Primary: 11S80. Secondary: 26A33.

\newpage

\section{Introduction}

{\bf 1.1.} In analysis of complex-valued functions on the field $\Qp$ of $p$-adic numbers (or, more generally, on a non-Archimedean local field), the basic operator is Vladimirov's fractional differentiation operator $\DA$, $\alpha >0$, defined via the Fourier transform or, for wider classes of functions, as a hypersingular integral operator \cite{K2001,VVZ}. Properties of this $p$-adic pseudo-differential operator were studied by Vladimirov (see \cite{VVZ}) and found to be more complicated than those of its classical counterparts. For example, as an operator on $L^2(\Qp )$, it has a point spectrum of infinite multiplicity. However, it was shown in \cite{K2014} to behave much simpler on radial functions $x\to f(|x|_p)$.

In particular, in \cite{K2014} the first author introduced a right inverse $I^\alpha$ to the operator $\DA$ on radial functions, which can be seen as a $p$-adic analog of the Riemann-Liouville fractional integral of real analysis (including the case $\alpha =1$ of the usual antiderivative). Just as the Riemann-Liouville fractional integral is a source of many problems of analysis, that must be true for the operator $I^\alpha$.

In this paper we study asymptotic properties of the function $I^\alpha f$ for a given asymptotic expansion of $f$; for the asymptotic properties of Riemann-Liouville fractional integral see \cite{R,SKM,W}.

\medskip
{\bf 1.2.} Let us recall the main definitions and notation used below.

Let $p$ be a prime
number. The field of $p$-adic numbers is the completion $\mathbb Q_p$ of the field $\mathbb Q$
of rational numbers, with respect to the absolute value $|x|_p$
defined by setting $|0|_p=0$,
$$
|x|_p=p^{-\nu }\ \mbox{if }x=p^\nu \frac{m}n,
$$
where $\nu ,m,n\in \mathbb Z$, and $m,n$ are prime to $p$. It is well known that $\Qp$ is a locally compact topological field with the topology determined by the metric $|x-y|_p$, and that there are no absolute values on $\mathbb Q$, which are not equivalent to the ``Euclidean'' one, or one of $|\cdot |_p$. We will denote by $dx$ the Haar measure on the additive group of $\mathbb Q_p$ normalized by the condition $\int_{|x|_p\le 1}dx=1$.

The absolute value $|x|_p$, $x\in \mathbb Q_p$, has the following properties:
\begin{gather*}
|x|_p=0\ \mbox{if and only if }x=0;\\
|xy|_p=|x|_p\cdot |y|_p;\\
|x+y|_p\le \max (|x|_p,|y|_p).
\end{gather*}
The latter property called the ultrametric inequality (or the non-Archi\-me\-dean property) implies the total disconnectedness of $\Qp$ and unusual geometric properties. Note also the following consequence of the ultrametric inequality:
\begin{equation*}
|x+y|_p=\max (|x|_p,|y|_p)\quad \mbox{if }|x|_p\ne |y|_p.
\end{equation*}

We will often use the integration formulas (see \cite{K2001,VVZ,V}):
$$
\int\limits_{|x|_p\le p^n}|x|_p^{\alpha -1}dx=\frac{1-p^{-1}}{1-p^{-\alpha }}p^{\alpha n};\text{ here and below }n\in \mathbb Z,\alpha >0;
$$
in particular,
$$
\int\limits_{|x|_p\le p^n}dx=p^n;
$$
$$
\int\limits_{|x|_p=p^n}dx=(1-\frac1p)p^n;
$$
$$
\int\limits_{|x|_p=1}|1-x|_p^{\alpha -1}=\frac{p-2+p^{-\alpha}}{p(1-p^{-\alpha})}.
$$

See \cite{K2001,VVZ} for further details of analysis of complex-valued functions on $\Qp$.

From now on, we consider the case $\alpha >1$. The integral operator $I^\alpha$ introduced in \cite{K2014} has the form
\begin{equation}
(I^\alpha f(x))=\frac{1-p^{-\alpha}}{1-p^{\alpha -1}}\int\limits_{|y|_p\le |x|_p}\left( |x-y|_p^{\alpha -1}-|y|_p^{\alpha -1}\right) f(y)\,dy,
\end{equation}
where $f$ is a locally integrable function on $\Qp$. See \cite{K2014} for its connection to the Vladimirov operator $\DA$ and applications to non-Archimedean counterparts of ordinary differential equations. Note that our results can be generalized easily to the case of general non-Archimedean local fields.

\medskip
\section{Asymptotics at the origin}

Let $0<M_0<M_1<M_2<\ldots$, $M_n\to \infty$. Then the sequence $f_n(x)=|x|_p^{M_n}$ is an asymptotic scale for $x\to 0$ (see, for example, \S 16 of \cite{SKM} for the main notions regarding asymptotic expansions).

\medskip
\begin{teo}
Suppose that a function $f$ admits an asymptotic series expansion
$$
f  \sim \sum\limits_{n=0}^\infty a_n |x|_p^{M_n}, \quad |x|_p \to 0, a_n\in \mathbb C.
$$
Then
\begin{equation}
(I^\alpha f(x))\sim \frac{1-p^{-\alpha}}{1-p^{\alpha -1}}\sum\limits_{n=0}^\infty a_nb_n |x|_p^{M_n+\alpha }, \quad |x|_p \to 0,
\end{equation}
where
$$
b_n=\frac{p^{-\alpha+1} -1}{(1-p^{-\alpha})p} +
(1-p^{-1})  \sum_{k=1}^\infty    (1-p^{-k(\alpha-1)} ) p^{-k (M_n+1)}.
$$
\end{teo}

\medskip
{\it Proof}. We have
\[
f  = \sum_{n=0}^N a_n |x|_p^{M_n} + R_N(x), \quad R_N(x) = o(|x|_p^{M_N}) \quad |x|_p \to 0.
\]
Then $I^\alpha f=I^\alpha_{(1)}+I^\alpha_{(2)}$,
$$
I^\alpha_{(1)}=\ufrac \int\limits_{|y|_p \leq |x|_p}
(|x-y|_p^{\alpha-1} - |y|_p^{\alpha-1}) \Bigl( \sum_{n=0}^N a_n |y|_p^{M_n} \Bigr)\,dy,
$$
$$
I^\alpha_{(2)}=\ufrac \int\limits_{|y|_p \leq |x|_p}
(|x-y|_p^{\alpha-1} - |y|_p^{\alpha-1})R_N(y)\,dy.
$$

After the change of variables $y=sx$ we get
$$
I^\alpha_{(1)}=\ufrac |x|_p^\alpha \int\limits_{|s|_p \leq 1}
(|1-s|_p^{\alpha-1} - |s|_p^{\alpha-1}) \Bigl( \sum_{n=0}^N a_n |x|_p^{M_n} |s|_p^{M_n} \Bigr) ds=\ufrac |x|_p^\alpha (A+B)
$$
where
\begin{multline*}
A=\int\limits_{|s|_p < 1}\left( 1-|s|_p^{\alpha-1} \right) \Bigl( \sum_{n=0}^N a_n |x|_p^{M_n} |s|_p^{M_n} \Bigr)\,ds\\
= \sum_{n=0}^N a_n |x|_p^{M_n} \sum_{k=1}^\infty  \left( 1-p^{-k(\alpha-1)} \right) p^{-k M_n} \int\limits_{|s|_p = p^{-k}}\,ds\\
= (1-p^{-1})\sum_{n=0}^N a_n |x|_p^{M_n}   \sum_{k=1}^\infty  \left( 1-p^{-k(\alpha-1)} \right) p^{-k (M_n+1)},
\end{multline*}
$$
B=\sum_{n=0}^N a_n |x|_p^{M_n} \int\limits_{|s|_p = 1}  \left( |1-s|_p^{\alpha-1} - 1\right)\,ds=\frac{p^{-\alpha+1} -1}{(1-p^{-\alpha})p}\sum_{n=0}^N a_n |x|_p^{M_n}.
$$

On the other hand, since $|R_N(x)|\le C|x|_p^{M_{N+1}}$, we find that for some constant $C_1>0$,
$$
\left| I^\alpha_{(2)}\right| \le C_1|x|_p^{\alpha+M_{N+1}}  \int\limits_{|s|_p \leq 1}  (|1-s|_p^{\alpha-1}-|s|_p^{\alpha-1})  |s|_p^{M_{N+1}} ds =
O(|x|_p^{\alpha+M_{N+1}}).
$$
The above calculations result in the asymptotic relation (2). $\qquad \blacksquare$

\medskip
\section{Asymptotics at infinity}

For positive functions $\varphi,\psi$, we write $\varphi (x)\asymp \psi (x)$, $|x|_p\to \infty$, if $c\psi (x)\le \varphi (x)\le d\psi (x)$, for large values of $|x|_p$, $x\in \Qp$, for some positive constants $c,d$.

\medskip
\begin{teo}
Suppose that $a\le f(x)\le b$ ($a,b>0$) for $|x|_p<1$, $|f(x)|\le C|x|_p^{-M}$, $M>1$, $C>0$, for $|x|_p\ge 1$. Then
\begin{equation}
(I^\alpha f)(x)\asymp |x|_p^{\alpha -1},\quad |x|_p\to \infty.
\end{equation}
\end{teo}

\medskip
{\it Proof}. Let us rewrite (1) with $|x|_p\ge 1$ in the form $I^\alpha f=J^\alpha_{(1)}f+J^\alpha_{(2)}f$ where
$$
\left( J^\alpha_{(1)}f\right) (x)=\frac{1-p^{-\alpha}}{1-p^{\alpha -1}}\int\limits_{|y|_p<1}\left( |x-y|_p^{\alpha -1}-|y|_p^{\alpha -1}\right) f(y)\,dy,
$$
$$
\left( J^\alpha_{(2)}f\right) (x)=\frac{1-p^{-\alpha}}{1-p^{\alpha -1}}\int\limits_{1\le |y|_p\le |x|_p}\left( |x-y|_p^{\alpha -1}-|y|_p^{\alpha -1}\right) f(y)\,dy
$$
Then
$$
\left( J^\alpha_{(1)}f\right) (x)\asymp \int\limits_{|y|_p<1}\left( |x-y|_p^{\alpha -1}-|y|_p^{\alpha -1}\right)\,dy\asymp |x|_p^{\alpha -1}.
$$

Next, if $|x|_p=p^N$, $N\ge 0$, then
\begin{multline*}
\left| \left( J^\alpha_{(2)}f\right) (x)\right| \le C\int\limits_{1\le |y|_p\le |x|_p}\left( |x-y|_p^{\alpha -1}-|y|_p^{\alpha -1}\right) |y|_p^{-M}dy\\
=C\left\{ \sum\limits_{j=0}^{N-1}\int\limits_{|y|_p=p^j}\left( |x|_p^{\alpha -1}-|y|_p^{\alpha -1}\right) |y|_p^{-M}dy
+\int\limits_{|y|_p=p^N}\left( |x-y|_p^{\alpha -1}-p^{N(\alpha -1)}\right) p^{-MN}dy\right\}\\
=C\left\{ (1-\frac1p )\sum\limits_{j=0}^{N-1}p^j\left( p^{N(\alpha -1)}-p^{j(\alpha -1)}\right) p^{-Mj}\right. \\
\left. +p^{-MN}\int\limits_{|y|_p=p^N}|x-y|_p^{\alpha -1}dy-(1-\frac1p ) p^{\alpha N-MN}\right\}.
\end{multline*}
Calculating the integral as above and finding the sums of geometric progressions we see that $\left| \left( J^\alpha_{(2)}f\right) (x)\right| \le \operatorname{const}\cdot |x|_p^{\alpha -1}$, which proves (3). $\qquad \blacksquare$

\medskip
\section{Logarithmic asymptotics}

If a function $f$ decays slower than it did under the assumptions of Theorem 2, then a richer asymptotic behavior is possible. Let us consider the case where $f(t)\ge 0$,
\begin{equation}
f(x)\sim |x|_p^{-\beta}\sum_{n=0}^\infty a_n  (\log |x|_p)^{\gamma-n} , \quad |x|_p \to \infty,
\end{equation}
where $0\le \beta <1$, $\gamma \ge 0$, $a_n\in \mathbb R$.

First we need some auxiliary results.

\medskip
\begin{lem}
Let $0\le f(x)=o\left( |x|_p^{-\lambda}\right)$, $|x|_p\to \infty$, where $0<\lambda <1$.
Then
\begin{equation}
G_1(r)\overset{\text{def}}{=}\int\limits_{|y|_p\le r}f(y)\,dy=o(r^{1-\lambda}),\quad r\to \infty.
\end{equation}
\end{lem}

\medskip
{\it Proof}. Let $n_0=[\log_p r]$. Then $p^{n_0}\le r\le p^{n_0+1}$. It is known (see Section 1) that
\begin{equation}
\int\limits_{|y|_p\le p^\nu }|y|_p^{-\lambda }dy=\frac{1-p^{-1}}{1-p^{\lambda -1 }}p^{(1-\lambda) \nu},\quad \nu \in \mathbb Z,
\end{equation}
so that
\begin{equation}
G_2(r)\overset{\text{def}}{=}\int\limits_{|y|_p\le r}|y|_p^{-\lambda}\,dy=O(r^{1-\lambda}),\quad r\to \infty.
\end{equation}

By our assumption, for any $n\in \mathbb N$, there exists such $r_0=r_0(n)$ that $f(x)<\frac1n |x|_p^{-\lambda}$ for $|x|_p>r_0$. Then we can write
$$
\frac{G_1(r)}{G_2(r)} =
\frac{G_1(r_0(n)) + ( G_1(r)-G_1(r_0(n)) )}{ G_2(r_0(n)) + ( G_2(r)-G_2(r_0(n)) ) }\le \frac{G_1(r_0(n)) + \frac1n G_3(n,r)}{ G_2(r_0(n)) + G_3(n,r)}
$$
where
$$
G_3(n,r)=\int\limits_{r_0\le |y|_p\le r}|y|_p^{-\lambda}\,dy.
$$

It follows from (6) that $G_3(n,r)\to \infty$, so that
$$
0\le \limsup\limits_{r\to \infty}\frac{G_1(r)}{G_2(r)}\le \frac1n
$$
where $n$ is arbitrary. Therefore
$$
\lim\limits_{r\to \infty}\frac{G_1(r)}{G_2(r)}=0,
$$
which gives, together with (7), the required asymptotic relation (5). $\qquad \blacksquare$

\bigskip
\begin{lem}
Let $0\le \beta <1$, $k\in \mathbb N$. For any $\varepsilon >0$, such that $\beta +\varepsilon <1$,
\begin{equation}
K_r\overset{\text{def}}{=}\int\limits_{|t|_p \leq r^{-1}}  ( |1-t|_p^{\alpha-1} - |t|_p^{\alpha-1} ) |t|_p^{-\beta} |\log |t|_p|^k dt =
O  \bigl(   r^{\beta+\varepsilon-1}\bigr),\quad r\to \infty.
\end{equation}
\end{lem}

\medskip
{\it Proof}. Assuming that $r>2$, we have $|t|_p<\frac12$, so that $|1-t|_p^{\alpha-1} - |t|_p^{\alpha-1}=1-|t|_p^{\alpha -1}\le 1$, and we find that
$$
K_r\le \int\limits_{|t|_p \leq r^{-1}}|t|_p^{-\beta} |\log |t|_p|^k dt\le \int\limits_{|t|_p \leq r^{-1}}|t|_p^{-\beta -\varepsilon }dt,
$$
if $r$ is large enough, and the relation (8) follows from the integration formula (6). $\qquad \blacksquare$

\bigskip
Now we are ready to consider the asymptotics of $I^\alpha f$ for a function $f$ satisfying (4). Below we use the notation
$$
\binom{\gamma}{n}=\frac{\gamma (\gamma -1)\cdots (\gamma -n+1)}{n!}
$$
for any real positive number $\gamma$ and $n\in \mathbb N$.

\medskip
\begin{teo}
If a function $f\ge 0$ satisfies the asymptotic relation (4), then
\begin{equation}
(I^{\alpha}f)(x) \sim  \ufrac |x|_p^{\alpha-\beta} \sum_{n=0}^\infty   B_n  ( \log |x|_p )^{\gamma-n}, \quad |x|_p \to \infty,
\end{equation}
where
$$
B_n =\sum_{k=0}^n   a_{n-k}\binom{\gamma+k-n}{k} \Omega(k,\alpha,\beta),
$$
$$
\Omega(k,\alpha,\beta)= \int\limits_{ |t|_p \leq 1} (|1-t|_p^{\alpha-1} - |t|_p^{\alpha-1})|t|_p^{-\beta}  (\log |t|_p)^{k} dt.
$$
\end{teo}

\medskip
{\it Proof}. Let us write $(I^{\alpha}f)(x)$ for $|x|_p\ge 1$ as the sum of two integrals $I_1$ and $I_2$, with the integration over $\{ y:\ |y|_p<|x|_p^{1/2}\}$ and $\{ y:\ |x|_p^{1/2}\le |y|_p\le |x|_p\}$ respectively.

Denote $\mathcal K(x,y)=|x-y|_p^{\alpha -1}-|y|_p^{\alpha -1}$. Considering $I_1$, for $|y|_p\le |x|_p$, we have
\begin{equation}
|\mathcal K(x,y)|\le |x|_p^{\alpha -1}.
\end{equation}
Indeed, if $|x|_p>1$, then $|y|_p<|x|_p$, $\mathcal K(x,y)=|x|_p^{\alpha -1}-|y|_p^{\alpha -1}$, and we get (10). If $|x|_p=1$, $|y|_p<1$, then $0<\mathcal K(x,y)=1-|y|_p^{\alpha -1}<|x|_p^{\alpha -1}$.

It follows from (10) that
$$
0\le I_1\le C|x|_p^{\alpha -1}\int\limits_{|y|_p<|x|_p^{1/2}}f(y)\,dy,
$$
and by (4) and Lemma 1, for any small $\varepsilon >0$,
\begin{equation}
I_1=o\left( |x|_p^{ \alpha - \beta +\frac{\beta+\varepsilon-1}{2} }\right),  \quad |x|_p \to \infty.
\end{equation}

Considering $I_2$ we write
$$
f(t) = |t|_p^{-\beta} \sum_{n=0}^N a_n  (\log |t|_p)^{\gamma-n} +R_N(t), \quad
R_N(t) = O(|t|_p^{-\beta} (\log |t|_p)^{\gamma-N-1}) \quad |t|_p \to \infty .
$$
Denote
\begin{multline*}
 L(\alpha,\beta,\gamma,x) = \int\limits_{|x|_p^{1/2} \leq |y|_p \leq |x|_p} \left( |x-y|_p^{\alpha-1} - |y|_p^{\alpha-1}\right) |y|_p^{-\beta}  (\log |y|_p)^{\gamma} \,dy \\
 =|x|_p^{\alpha-\beta}  ( \log |x|_p )^{\gamma} \int\limits_{|x|_p^{-1/2} \leq |t|_p \leq 1} \left( |1-t|_p^{\alpha-1} - |t|_p^{\alpha-1}\right) |t|_p^{-\beta} \Bigl( 1 + \frac{\log |t|_p}{\log |x|_p} \Bigr) ^{\gamma} dt
\end{multline*}
where on the domain of integration,
$$
\left| \frac{\log |t|_p}{\log |x|_p}\right| \le \frac12,
$$
and we may write, for a non-integer $\gamma$, the convergent binomial series
$$
 \Bigl( 1 + \frac{\log |t|_p}{\log |x|_p} \Bigr) ^{\gamma} =\sum\limits_{k=0}^\infty \binom{\gamma}{k}\left( \frac{\log |t|_p}{\log |x|_p}\right)^k.
$$

Note that we can use the Taylor formula with the integral form of the remainder
$$
(1+s)^\gamma =\sum\limits_{k=0}^N \binom{\gamma}{k}s^k+\frac{\gamma (\gamma -1)\cdots (\gamma -N)}{N!}\int\limits_0^s (1+\sigma )^{\gamma -N-1}(s-\sigma )^N\,d\sigma
$$
where
$$
\int\limits_0^s (1+\sigma )^{\gamma -N-1}(s-\sigma )^N\,d\sigma =s^{N+1}\int\limits_0^1 (1+s\tau )^{\gamma -N-1}(1-\tau )^N\,d\tau =s^{N+1}\int\limits_0^1 (1+s(1-\tau ))^{\gamma -N-1}\tau^N\,d\tau.
$$
If $-\frac12 <s<\frac12$, $0<\tau <1$, then $\frac12 \le 1+s(1-\tau )\le \frac32$. Therefore
$$
 \Bigl( 1 + \frac{\log |t|_p}{\log |x|_p} \Bigr) ^{\gamma} =\sum\limits_{k=0}^N\binom{\gamma}{k}\left( \frac{\log |t|_p}{\log |x|_p}\right)^k+S_N(t,x),
$$
$$
S_N(t,x)=O\left( \left( \frac{\log |t|_p}{\log |x|_p}\right)^{N+1}\right) ,\quad |x|_p\to \infty,
$$
and this asymptotics is uniform with respect to $t,|t|_p\in [|x|_p^{-1/2},1]$.

Substituting and using Lemma 2 we obtain the expansion
\begin{multline}
L(\alpha,\beta,\gamma,x) = |x|_p^{\alpha-\beta} \sum_{k=0}^N  \binom{\gamma}{k}  \Omega(k,\alpha,\beta)( \log |x|_p )^{\gamma-k}\\
+o(|x|_p^{\alpha-\beta}(\log |x|_p)^{\gamma-N}), \quad |x|_p \to \infty.
\end{multline}

We have
$$
 I_2=\ufrac \sum_{n=0}^N a_n L(\alpha,\beta,\gamma-n,x) +\ufrac
\int\limits_{|x|_p^{1/2} \leq |y|_p \leq |x|_p} \left( |x-y|_p^{\alpha-1} - |y|_p^{\alpha-1}\right) R_N(y)\,dy
$$
where
\begin{multline*}
\int\limits_{|x|_p^{1/2} \leq |y|_p \leq |x|_p} \left( |x-y|_p^{\alpha-1} - |y|_p^{\alpha-1}\right) R_N(y)\,dy \\
\le C L(\alpha ,\beta ,\gamma -N-1,x)=O\left( |x|_p^{\alpha -\beta}\left( \log |x|_p\right)^{\gamma -N-1}\right),\quad |x|_p\to \infty.
\end{multline*}
The last estimate is a consequence of (12).

Now the asymptotic relations (11) and (12) imply the required relation (9). $\qquad \blacksquare$

\medskip
In our final result, we give a modification of Theorem 3 for the case where $\beta =1$.

\medskip
\begin{teo}
Suppose that $f$ is nonnegative,
$$
f(x) \sim |x|_p^{-1} \sum_{n=0}^\infty  a_n   (\log |x|_p)^{\gamma-n} ,\quad |x|_p \to \infty.
$$
Then
\begin{equation}
(\iafi)(x) \sim \ufrac \Bigl[  |x|_p^{\alpha-1} \int\limits_{|y|_p \leq |x|_p} f(y) dy+\sum_{n=0}^\infty \widetilde{B}_n( \log |x|_p )^{\gamma-n}
\Bigr], \quad |x|_p \to \infty
\end{equation}
where
$$
\widetilde{B}_n=\sum_{k=0}^n  a_{n-k} \binom{\gamma+k-n}{k} \widetilde{\Omega}(k,\alpha),
$$
$$
\widetilde{\Omega}(k,\alpha)=\int\limits_{ |t|_p \leq 1}   \left( |1-t|_p^{\alpha-1} - |t|_p^{\alpha-1}-1\right) |t|_p^{-1}  (\log |t|_p)^k\, dt.
$$
\end{teo}

\medskip
{\it Proof}. Let us write $I^\alpha f= \ufrac (J_1+J_2+J_3)$ where
$$
J_1=\int\limits_{|y|_p \leq |x|_p^{1/2}} \left( |x-y|_p^{\alpha-1} - |y|_p^{\alpha-1} - |x|_p^{\alpha-1}\right) f(y)\,dy,
$$
$$
J_2=\int\limits_{|x|_p^{1/2} \leq |y|_p \leq |x|_p} \left( |x-y|_p^{\alpha-1} - |y|_p^{\alpha-1} - |x|_p^{\alpha-1}\right) f(y)\,dy,
$$
$$
J_3=|x|_p^{\alpha-1}\int\limits_{|y|_p \leq |x|_p}f(y)\,dy .
$$

Choosing $\varepsilon >0$, such that $1+\varepsilon <\alpha$, we see that $f(x)=o\left( |x|_p^{-1+\varepsilon}\right)$, $|x|_p\to \infty$. By
Lemma 1,
$$
\int\limits_{|y|_p \leq |x|_p^{1/2}}f(y) dy  =o \bigl(  |x|_p^{\frac{\varepsilon}{2}}) \bigr) , \quad |x|_p \to \infty.
$$

For the kernel of the above integral operator we get, considering various cases, the estimate
$$
| |x-y|_p^{\alpha-1} - |y|_p^{\alpha-1} - |x|_p^{\alpha-1} |  \leq 2 |y|_p^{\alpha-1}
$$
It follows from Lemma 1 that
\begin{equation}
|J_1|\le 2\int\limits_{|y|_p \leq |x|_p^{1/2}} |y|_p^{\alpha-1}f(y)\, dy = o \bigl( |x|_p^{\frac{\alpha-1+\varepsilon}{2}}\bigr), \quad |x|_p \to \infty.
\end{equation}

By our assumption,
$$
f(t) =|t|_p^{-1} \sum_{n=0}^N  a_n   (\log |t|_p)^{\gamma-n}
+ R_N(t),\quad
R_N(t) = O \bigl( |t|_p^{-1} (\log |t|_p)^{\gamma-N-1} \bigr) ,\quad |t|_p \to \infty .
$$
Let us consider the expression
\begin{multline*}
\widetilde{L} (\alpha,\gamma,x) =
\int\limits_{|x|_p^{1/2} \leq |y|_p \leq |x|_p} \left( |x-y|_p^{\alpha-1} - |y|_p^{\alpha-1} - |x|_p^{\alpha-1}\right) |y|_p^{-1} (\log|y|_p)^{\gamma}  dy \\
=|x|_p^{\alpha-1} \int\limits_{|x|_p^{-1/2} \leq |t|_p \leq 1}   (|1-t|_p^{\alpha-1} - |t|_p^{\alpha-1} - 1) |t|_p^{-1}  (\log |x|_p + \log |t|_p)^{\gamma}\, dt.
\end{multline*}

It follows from the first integration formula from Section 1 that
$$
\int\limits_{|t|_p \leq |x|_p^{-1/2} } (|1-t|_p^{\alpha-1} - |t|_p^{\alpha-1} - 1) |t|_p^{-1} (\log |t|_p)^{k} \,dt =o \bigl( |x|_p^{\frac{1-\alpha+\varepsilon}{2}}
\bigr), \quad |x|_p \to \infty.
$$
This implies (just as in the proof of Theorem 3) the expansion
$$
\widetilde{L}(\alpha,\gamma,x) \sim |x|_p^{\alpha-1} \sum_{k=0}^\infty \binom{\gamma}{k} ( \log |x|_p )^{\gamma-k} \widetilde{\Omega}(k,\alpha) , \quad
|x|_p \to \infty.
$$
Taking into account (14), we come to (13). $\qquad \blacksquare$

\medskip
\section*{Acknowledgments}
The work of the first author was supported in part by Grant 23/16-18 ``Statistical dynamics, generalized Fokker-Planck equations, and their applications in the theory of complex systems'' of the Ministry of Education and Science of Ukraine.


\medskip
\begin{thebibliography}{999}

\bibitem{K2001}
A. N. Kochubei, {\it Pseudo-Differential Equations and Stochastics
over Non-Archimedean Fields}, Marcel Dekker, New York, 2001.
\bibitem{K2014}
A. N. Kochubei, Radial solutions of non-Archimedean pseudo-differential equations, {\it Pacif. J. Math.} {\bf 269} (2014), 355--369.
\bibitem{R}
E. Rieksti\d{n}\u{s}, Asymptotic representation of certain types of the convolution integral, {\it Latviiski Matem. Ezhegodnik} {\bf 8} (1970), 223--239 (Russian).
\bibitem{SKM}
S. G. Samko, A. A. Kilbas, and O. I. Marichev, {\it Fractional Integrals and Derivatives: Theory and Applications}, Gordon and Breach, New York, 1993.

\bibitem{VVZ}
V. S. Vladimirov, I. V. Volovich and E. I. Zelenov, {\it $p$-Adic Analysis and
Mathematical Physics}, World Scientific, Singapore, 1994.
\bibitem{V}
V. S. Vladimirov, {\it Tables of Integrals of Complex-Valued
Functions of $p$-Adic Arguments}, Steklov Mathematical Institute,
Moscow, 2003 (Russian). English version, ArXiv: math-ph/9911027.
\bibitem{W}
R. Wong, Asymptotic expansions of fractional integrals involving logarithms, {\it SIAM J. Math. Anal.} {\bf 9} (1978), 835--842.
\end{thebibliography}
\end{document}